\documentclass[12pt]{amsart}
 \usepackage{amscd,amsmath,amsthm,amssymb}
 \usepackage{pstcol,pst-plot,pst-3d}
 \psset{unit=0.7cm,linewidth=0.8pt,arrowsize=2.5pt 4}

\newpsstyle{fatline}{linewidth=1.5pt}
\newpsstyle{fyp}{fillstyle=solid,fillcolor=verylight}
\definecolor{verylight}{gray}{0.97}
\definecolor{light}{gray}{0.93}
\definecolor{medium}{gray}{0.82}
 \def\frk{\frak}               

 \def\mm{{\frk m}}

 \def\Sc{{\mathcal S}}

 \def\xb{{\bold x}}

 \def\opn#1#2{\def#1{\operatorname{#2}}} 
 \opn\chara{char} \opn\length{\ell} \opn\pd{pd} \opn\rk{rk}
 \opn\projdim{proj\,dim} \opn\injdim{inj\,dim} \opn\rank{rank}
 \opn\depth{depth} \opn\grade{grade} \opn\height{height}
 \opn\embdim{emb\,dim} \opn\codim{codim}
 
 \opn\Tr{Tr} \opn\bigrank{big\,rank}
 \opn\superheight{superheight}\opn\lcm{lcm}
 \opn\trdeg{tr\,deg}
 \opn\reg{reg} \opn\lreg{lreg} \opn\ini{in} \opn\lpd{lpd}
 \opn\size{size} \opn\sdepth{sdepth}
 \opn\link{link}\opn\fdepth{fdepth}\opn\lex{lex}
 \opn\div{div} \opn\Div{Div} \opn\cl{cl} \opn\Cl{Cl}
 \opn\Spec{Spec} \opn\Supp{Supp} \opn\supp{supp} \opn\Sing{Sing}
 \opn\Ass{Ass} \opn\Min{Min}\opn\Mon{Mon}
 \opn\Ann{Ann} \opn\Rad{Rad} \opn\Soc{Soc}
 \opn\Im{Im} \opn\Ker{Ker} \opn\Coker{Coker} \opn\Am{Am}
 \opn\Hom{Hom} \opn\Tor{Tor} \opn\Ext{Ext} \opn\End{End}
 \opn\Aut{Aut} \opn\id{id}
 
 \opn\nat{nat}
 \opn\pff{pf}
 \opn\Pf{Pf} \opn\GL{GL} \opn\SL{SL} \opn\mod{mod} \opn\ord{ord}
 \opn\Gin{Gin} \opn\Hilb{Hilb}\opn\sort{sort}
 \opn\aff{aff} \opn
 \con{conv} \opn\relint{relint} \opn\st{st}
 \opn\lk{lk} \opn\cn{cn} \opn\core{core} \opn\vol{vol}
 \opn\link{link} \opn\star{star}\opn\lex{lex}\opn\set{set}
 \opn\gr{gr}
 
 \def\pot#1#2{#1[\kern-0.28ex[#2]\kern-0.28ex]}

 \opn\dirlim{\underrightarrow{\lim}}
 \opn\inivlim{\underleftarrow{\lim}}
 \let\sect=\cap
 \let\dirsum=\oplus

 \let\Dirsum=\bigoplus
 
 \let\to=\rightarrow
 
 \def\Implies{\ifmmode\Longrightarrow \else
         \unskip${}\Longrightarrow{}$\ignorespaces\fi}
 \def\implies{\ifmmode\Rightarrow \else
         \unskip${}\Rightarrow{}$\ignorespaces\fi}
 \def\iff{\ifmmode\Longleftrightarrow \else
         \unskip${}\Longleftrightarrow{}$\ignorespaces\fi}

 \let\:=\colon
 \newtheorem{Theorem}{Theorem}

 \let\epsilon\varepsilon
 \let\kappa=\varkappa
 \def\qed{\ifhmode\textqed\fi
       \ifmmode\ifinner\quad\qedsymbol\else\dispqed\fi\fi}
 \def\textqed{\unskip\nobreak\penalty50
        \hskip2em\hbox{}\nobreak\hfil\qedsymbol
        \parfillskip=0pt \finalhyphendemerits=0}
 \def\dispqed{\rlap{\qquad\qedsymbol}}
 
 \opn\dis{dis}
 \def\pnt{{\raise0.5mm\hbox{\large\bf.}}}
 
 \opn\Lex{Lex}

\begin{document}

 \title {Powers are Golod}

 \author {J\"urgen Herzog}

\address{J\"urgen Herzog, Fachbereich Mathematik, Universit\"at Duisburg-Essen, Campus Essen, 45117
Essen, Germany}
\email{juergen.herzog@uni-essen.de}

\subjclass[2000]{13A02, 13D40}
\keywords{Powers of ideals, Golod rings, Koszul cycles}

 \begin{abstract}
Let $I$ be a proper graded ideal in a positively  graded polynomial ring $S$ over a field of characteristic $0$. In this note it is shown that $S/I^k$ is Golod for all $k\geq 2$.
 \end{abstract}

\thanks{The paper was written while the author was visiting MSRI at Berkeley. He thanks for the support, the hospitality and the inspiring atmosphere of this institution}

 \maketitle

 \section*{Introduction}
 Let $(R,\mm)$ be a Noetherian local  ring with residue class field $K$, or a standard graded $K$-algebra with graded maximal ideal $\mm$. The formal power series $P_R(t)= \sum_{i \geq 0} \dim_K \Tor_i^{R} (R/\mm,R/\mm) t^i$ is called the {\em Poincar\'{e}  series} of $R$. Though the ring is Noetherian, the Poincar\'{e}  series of $R$  is in general not a rational function. The first example that showed that $P_R(t)$ may  not rational was given by Anik \cite{An}.  In the meantime more such examples are known, see \cite{Ro} and its references. On the other hand, Serre showed that $P_R(t)$ is coefficientwise bounded above by the rational series
 \[
 \frac{(1+t)^n}{1-t\sum_{i\geq 0}\dim_K H_i(\xb;R)t^i},
 \]
where $\xb=x_1,\ldots,x_n$ is a minimal system of generators of $\mm$ and where $H_i(\xb;R)$ denotes the $i$th Koszul homology of the sequence $\xb$.

The ring $R$ is called {\em Golod}, if $P_R(t)$ coincides with this upper bound given by Serre. There is also a relative version of Golodness which is defined for local homomorphisms as an obvious extension of the above concept of Golod rings. We refer the reader for details regarding Golod rings and Golod homomorphism to the survey article \cite{Av} by Avramov. Here we just want to quote the following result of Levin \cite{Le} which says that for any Noetherian local ring $(R,\mm)$,  the canonical map $R\to R/\mm^k$ is a Golod homomorphism for all $k\gg 0$. It is natural to ask whether in this statement $\mm$ could be replaced by any other proper ideal of $I$. Some very recent results indicate that this question may have a positive answer. In fact, in \cite{HWY} it is shown that if $R$ is regular, then for any proper ideal $I\subset R$ the residue class ring  $R/I^k$ is Golod for $k\gg 0$, which, since $R$ is regular, is equivalent to  saying that the residue class map $R\to R/I^k$ is a Golod homomorphism  for $k\gg 0$. But how big $k$ has to be chosen to make sure that $R/I^k$ is Golod? In the case that $R$ is the polynomial ring and $I$ is a proper monomial ideal, the surprising answer is that $R/I^k$ is Golod for all $k\geq 2$, as has been shown by Fakhari and Welker  in \cite{FW}. The authors show even more: if $I$ and $J$ a proper monomial ideals, then $R/IJ$ is Golod. Computational evidence as well as  a result of Huneke \cite{Hu} which says that for a regular local ring $R$, the residue class ring $R/IJ$ is never Gorenstein, unless $I$ and $J$ are principal ideals, suggest that $R/IJ$ is Golod  for any two proper ideals $I,J$ in a local ring (or graded ideals in  graded ring). Indeed, being Golod implies in particular that the Koszul homology $H(\xb;R)$ admits trivial multiplication,  while for a Gorenstein ring, by a result of Avramov and Golod \cite{LAv},  the multiplication map induces for all $i$ a non-degenerate pairing $H_i(\xb;R)\times H_{p-i}(\xb;R) \to H_p(\xb:R)$ where $p$ is the top non-vanishing homology of the Koszul homology. In the case that $I$ and $J$ are not necessarily monomial ideals, it is only known that $R/IJ$ is Golod if $IJ=I\sect J$, see \cite{HSt}.

\medskip
In the present note  we prove

\begin{Theorem}
\label{one}
Let $K$   be a field of characteristic $0$ and $S=K[x_1,\ldots,x_n]$ the graded polynomial ring over $K$ with $\deg x_i=a_i>0$ for $i=1,\ldots,n$. Let $I\subset S$ be a  graded ideal different from $S$. Then $S/I^k$ is Golod for all $k\geq 2$.
\end{Theorem}

It should be noted that this theorem does not imply the result of Fakhari and Welker, since here we consider only powers and not products of ideals,  and moreover we have to require that the base field is of characteristic $0$. Actually as the proof will show, it is enough to require in Theorem~\ref{one} that the characteristic of $K$ is big enough compared with the shifts in the graded free resolution of the ideal.

\section{Proof of Theorem 1}

Let $J\subset S_K[x_1,\ldots,x_n]$ be a graded ideal different from $S$ and set $R=S/J$. We denote by $(K(R), \partial)$ the Koszul complex of $R$ with respect to the sequence  $\xb=x_1,\ldots,x_n$. Let  $Z(R)$, $B(R)$ and $H(R)$ denote the module of cycles, boundaries and the homology of $K(R)$.

Golod \cite{Go} showed that Serre's upper bound for the Poincar\'{e} series
 is reached if and only if all Massey operations of $R$ vanish. By definition, this is the case (see \cite[Def. 5.5 and 5.6]{AKM}), if for each subset $\mathcal{S}$ of
  homogeneous elements of $\Dirsum_{i=1}^nH_i(R)$  there exists a function $\gamma$, which is defined on the set of finite
  sequences of elements from $\Sc$  with values in $\mm\dirsum\Dirsum_{i=1}^nK_i(R)$,  subject to the following conditions:
  \begin{enumerate}
    \item[(G1)] if $h\in \Sc$, then $\gamma(h)\in Z(R)$ and $h=[\gamma(h)]$;
    \item[(G2)] if  $h_1,\ldots,h_m$ is a sequence in $\mathcal{S}$ with $m>1$,  then
    \[
      \partial\gamma(h_1,\ldots,h_m)=\sum_{\ell=1}^{m-1}\overline{\gamma(h_1,\ldots,h_\ell)}\gamma(h_{\ell+1},\ldots,h_m),
    \]
    where $\bar{a} = (-1)^{i+1}a$ for $a\in K_i(R)$.
  \end{enumerate}

Note that (G2) implies,  that $\gamma(h_1)\gamma(h_2)$ is a boundary for all $h_1,h_2\in \Sc$ (which in particular implies the Koszul homology of Golod ring has trivial multiplication).  Suppose now that for each $\Sc$ we can choose a functions $\gamma$ such that $\gamma(h_1)\gamma(h_2)$ is not only  a boundary  but that $\gamma(h_1)\gamma(h_2)=0$ for all $h_1,h_2\in \Sc$.  Then obviously we may set $\gamma(h_1,\ldots h_r)=0$  for all $r\geq 2$, so that in this case (G2) is satisfied and $R$ is Golod.

The proof of Theorem~\ref{one} follows, once we have shown that $\gamma$ can be chosen that $\gamma(h_1)\gamma(h_2)=0$ for all $h_1,h_2\in \Sc$, in the case that  $J=I^k$ where $I$ is a graded ideal and $k\geq 2$. For the proof of this fact we use the following result from \cite{H}: Let
\[
0\to F_p\to F_{n-1}\to \cdots \to F_1\to F_0\to S/J
\]
be the graded minimal free $S$-resolution of $S/J$, and for each $i$ let $f_{11},\ldots, f_{ib_i}$ a homogeneous basis of $F_i$. Let $\varphi_i\: F_i\to F_{i-1}$ denote the chain maps in the resolution, and let
\[
\varphi_i(f_{ij})=\sum_{k=1}^{b_{i-1}}\alpha_{jk}^{(i)}f_{i-1,k},
\]
where the $\alpha_{jk}^{(i)}$ are  homogeneous polynomials.

In \cite[Corollary 2]{H} it is shown that for all $l=1,\ldots,p$  the elements
\[
\sum_{1\leq i_1<i_2<\cdots <i_l\leq n}a_{i_1}a_{i_2}\cdots a_{i_l}\sum_{j_2=1}^{b_{l-1}}\cdots \sum_{j_l=1}^{b_1}c_{j_1,\ldots,j_l}\frac{\partial(\alpha_{j_1,j_2}^{(l)},\alpha_{j_2,j_3}^{(l-1)},\ldots, \alpha_{j_l,1}^{(1)})}{\partial(x_{i_1},\ldots,x_{i_l})}e_{i_1}\wedge \cdots \wedge e_{i_l},
\]
$j_1=1,\ldots,b_l$ are cycles of $K(R)$ whose homology classes form a $K$-basis of $H_l(R)$.

Thus we see that a $K$-basis of $H_l(R)$ is given by cycles which are linear combinations of Jacobians determined by the entries $\alpha_{jk}^{(i)}$  of the matrices describing the resolution of $S/J$. The coefficients $c_{j_1,\ldots,j_l}$ which appear in these formulas are rational numbers determined by the degrees of the $\alpha_{jk}^{(i)}$, and the elements $e_{i_1}\wedge \cdots \wedge e_{i_l}$ form the $R$-basis of $K_l(R)=\bigwedge^l(\Dirsum_{i=1}^nRe_i)$.

From this result it follows that any homology class of $H_l(R)$ can be represented by a cycle which is a linear combination of Jacobians of the form

\begin{eqnarray}
\label{jacobian}
\frac{\partial(\alpha_{j_1,j_2}^{(l)}\alpha_{j_2,j_3}^{(l-1)},\ldots, \alpha_{j_l,1}^{(1)})}{\partial(x_{i_1},\ldots,x_{i_l})}.
\end{eqnarray}
We choose such representatives  for the elements of the set $\Sc$. Thus we may choose the map $\gamma$ in such a way that it assigns to each element of $\Sc$ a cycle which is a linear combination of Jacobians as in  (\ref{jacobian}).

Let $I$ by minimally generated by the polynomials $g_1,\ldots,g_r$. Then for $J=I^k$ we may choose a minimal set of generators  whose elements are of the form $g_{j_1}g_{j_2}\cdots g_{j_k}$. Thus in the Jacobians given in (1), each of the  entries $\alpha_{j_l,1}^{(1)}$ is a product of $k$ generators of $I$.

Since
\begin{eqnarray*}
\frac{\partial(*,\ldots, *, g_{j_1}g_{j_2}\cdots g_{j_k})}{\partial(x_{i_1},\ldots,x_{i_l})}=\sum_{s=1}^k\frac{\partial(*,\ldots, *, g_s)}{\partial(x_{i_1},\ldots,x_{i_l})}g_{j_1}\cdots g_{j_{s-1}} g_{j_{s+1}}\cdots g_{j_k},
\end{eqnarray*}
it follows that $\gamma(h)\in Z(R)\sect I^{k-1}K(R)$. From this it follows immediately that $\gamma(h_1)\gamma(h_2)=0$, as desired.

\end{document}